\documentclass[12pt]{article}

\usepackage{latexsym}
\usepackage{calc}
\usepackage{datetime}
\usepackage{amssymb, amsmath, amsfonts, listings, mathrsfs} 
\usepackage{xcolor} 
\usepackage[normalem]{ulem} 

\usepackage{graphicx}
\usepackage[labelfont=bf]{caption}

\newcounter{hours}\newcounter{minutes}

\newcommand{\stkout}[1]{\ifmmode\text{\sout{\ensuremath{#1}}}\else\sout{#1}\fi} 

\def\nr{\par \noindent}

\def\Def{\stackrel{\mathrm{def}}{=}}

\def\dom{{\rm dom \,}}
\def\vf{\varphi}
\def\beq{\begin{equation}}
\def\eeq{\end{equation}}

\def\R{\mathbb{R}}
\def\E{\mathbb{E}}

\def\BI{\begin{itemize}}
\def\EI{\end{itemize}}
\def\II{\item}
\newcommand{\SetEQ}{\setcounter{equation}{0}}
\newcommand{\refLE}[1]{\ensuremath{\stackrel{(\ref{#1})}{\leq}}}
\newcommand{\refEQ}[1]{\ensuremath{\stackrel{(\ref{#1})}{=}}}
\newcommand{\refGE}[1]{\ensuremath{\stackrel{(\ref{#1})}{\geq}}}

\newtheorem{theorem}{Theorem}
\newtheorem{lemma}{Lemma}
\newtheorem{corollary}{Corollary}

\newtheorem{assumption}{Assumption}
\newtheorem{definition}{Definition}

\newtheorem{example}{Example}
\newtheorem{remark}{Remark}
\newcommand{\proof}{\bf Proof: \rm \nr}
\newcommand{\qed}{\hfill $\Box$ \nr \medskip}
\newcommand{\half}{\mbox{$\frac{1}{2}$}}

\def\ba{\begin{array}}
\def\ea{\end{array}}
\def\beann{\begin{eqnarray*}}
\def\eeann{\end{eqnarray*}}
\def\bea{\begin{eqnarray}}
\def\eea{\end{eqnarray}}

\def\BT{\begin{theorem}}
\def\ET{\end{theorem}}
\def\BL{\begin{lemma}}
\def\EL{\end{lemma}}
\def\BC{\begin{corollary}}
\def\EC{\end{corollary}}
\def\BE{\begin{example}}
\def\EE{\end{example}}
\def\BD{\begin{definition}}
\def\ED{\end{definition}}
\def\BR{\begin{remark}}
\def\ER{\end{remark}}
\def\BAS{\begin{assumption}}
\def\EAS{\end{assumption}}
\def\BI{\begin{itemize}}
\def\EI{\end{itemize}}

\def\BMP{\begin{minipage}{9.5cm}}
\def\EMP{\end{minipage}}
\def\MPT{\begin{minipage}{11.5cm}}
\def\EPT{\end{minipage}}

\def\la{\langle}
\def\ra{\rangle}

\def\QF{\hspace{5ex} \Box}
\def\QR{\hfill \Box}

\title{
{\normalsize CORE DISCUSSION PAPER }\\{\normalsize
2022/01}\\\vspace{10mm} \textbf{Quartic Regularity 
}}
\author{Yu. Nesterov
\thanks{Center for Operations Research and Econometrics (CORE),
Catholic University of Louvain (UCL). E-mail:
Yurii.Neterov@uclouvain.be. This paper has received funding from the European Research Council (ERC) under the European Union’s
Horizon 2020 research and innovation programme (grant agreement No 788368).
It was also funded by Multidisciplinary Institute in Artificial intelligence MIAI@Grenoble Alpes (ANR-19-P3IA-0003).
}}

\date{
January 13, 2022
}

\begin{document}
\maketitle

\abstract{In this paper, we propose new linearly convergent second-order methods for minimizing convex quartic polynomials. This
framework is applied for designing optimization schemes, which can solve general convex problems satisfying a new condition of quartic regularity. It assumes positive definiteness and boundedness of the fourth derivative of the objective function. For such problems, an appropriate quartic regularization of Damped Newton Method has global linear rate of convergence. We discuss several important consequences of this result. In particular, it can be used for constructing new second-order methods in the framework of high-order proximal-point schemes \cite{High,Prox2}. These methods have convergence rate $\tilde O(k^{-p})$, where $k$ is the iteration counter, $p$ is equal to 3, 4, or 5, and tilde indicates the presence of logarithmic factors in the complexity bounds for the auxiliary problems, which are solved at each iteration of the schemes. }

\vspace{10ex}\noindent

\thispagestyle{empty}

\newpage\setcounter{page}{1}

\section{Introduction}
\setcounter{equation}{0}

\vspace{1ex}\noindent
{\bf Motivation.} After appearance in 2006 of the first second-order method with global efficiency bounds \cite{CNM}, the question of finding the right framework for complexity analysis of the higher-order schemes became one of the most interesting research directions in Convex Optimization. Indeed, the performance estimates of Cubic Newton Method \cite{CNM} and its accelerated version \cite{Accel} clearly demonstrated that the existing classification of complexity of optimization problems in accordance to the {\em condition number} of the objective\footnote{$^)$ It is traditionally defined by the uniform lower and upper bounds for the eigenvalues of the Hessian of the objective function. These constants control the growth of the objective with respect to its linear approximation at the current test point.}$^)$ fits well only the first-order schemes. The higher-order methods are insensitive to this characteristics since, in some sense, the influence of the quadratic part of the objective is eliminated by a single iteration of the simplest second-order scheme. Thus, it became necessary to find new characteristics helping in ranking the performance of the new methods.

The first attempt of introducing a notion of the second-order non-degeneracy was done already in \cite{Accel}. In the complexity analysis of \cite{CNM}, the assumption for bounding the growth of the objective was taken from the classical result of Kantorovich \cite{Kant} on local quadratic convergence of Newton Method. This was the Lipschitz continuity of the Hessian of the objective function. However, in \cite{Accel}, the lower bound for the growth of the objective was ensured by uniform convexity of {\em degree three}. The justification of this combination was purely algebraic: under these assumptions it was possible to prove the global {\em linear} rate of converge of the second-order methods. As compared to the theory of first-order methods, this choice looks much less adequate since it is not supported by any regular source of such problems, related to the real-life applications. 

Still, during the next decade, the development of the second-order methods was accomplished within this framework. For completeness of the picture, we mention interesting results on different variants of Cubic Regularization \cite{CC1,CC2}, universal methods \cite{GG1}, accelerated second-order methods based on contractions \cite{ND1,ND2}, and many others.

However, very soon it was shown that for convex functions we can efficiently implement also the higher-order optimization methods \cite{ImpTensor}. This result justifies an extensive study of these schemes \cite{YN1,YN2,Bubeck,Gas}, including the lower complexity bounds \cite{LB,LB2}. In particular, it was shown that the methods of the order three and higher have much faster rate of convergence than the best possible rate of the second-order schemes (on their ``natural'' problem classes).

The latter fact explains the level of the interest caused by an observation that practically all third-order methods can be implemented by a proper use of the second-order oracle, preserving at the same time their higher rate of convergence \cite{Super, High,Prox2}. The first results of this type were based on the difference approximation of the product of the third-order derivative by two vectors \cite{Super}, while the other papers employ the framework of the high-order proximal-point operators.

These results reveal the question on a proper attribution of optimization methods to particular problem classes. This link can be justified only by comparing the rates of convergence of the optimization schemes with the lower complexity bounds. Today it is clear that the initial naive hope that every problem class (formed by functions with uniformly bounded derivative of certain degree) can be unanimously related to the methods of certain order is too optimistic. This is may be true for the first-order methods, traditionally assigned to functions with bounded second derivatives. However, the next step in this classification is already much more involved.

Indeed, if we agree that the natural problem class for the second-order methods is formed by functions with bounded third derivative, then we must accept the corresponding bound on the maximal possible global rate of convergence in functional residual, which is of the order $O(k^{-3.5})$, with $k$ being the number of calls of oracle \cite{LB2}. On the other hand, as it was shown in \cite{ImpTensor}, the third derivative of convex functions is {\em not} an independent characteristic. It can be estimated from above by a combination of the second and the fourth derivatives. And if we assume instead the boundedness of these two derivatives, then the limits in the rate of convergence are relaxed up to $O(k^{-5})$. 

From the view point of global efficiency bounds, this perspective is very interesting. One possible way of implementing this idea was already presented in the papers \cite{High,Prox2}. The corresponding approach is based on the notion of the high-order proximal operators, approximated by auxiliary minimization schemes guided by the relative smoothness condition \cite{BBT,LFN}. In this paper, we develop a more direct approach based on the bounds for the {\em fourth derivative}. Our main observation is that the fourth derivative (contrary to the third one) can be {\em positive semidefinite}, and this property helps in further acceleration of optimization schemes.

\vspace{1ex}\noindent
{\bf Contents.} In Section \ref{sc-CQP}, we consider the problem of minimizing convex quartic polynomial in composite setting. We show that the corresponding quartic form is always convex and study its properties. In the next Section \ref{sc-Min}, we present the Damped Quartic Newton Method (DQNM), which regularizes the quadratic approximation of the objective by the fourth power of the norm, defined by the forth derivative. We prove that the corresponding global rate of convergence is linear with the absolute reduction factor. In Section \ref{sc-Simple}, we show that the regularization norm in DQNM can be replaced by a simpler one. The rate of convergence of the corresponding scheme is still linear, but it depends on the condition number of the new norm with respect to the optimal one. Note that the corresponding variant of DQNM is essentially a second-order method.

In Section \ref{sc-QReg}, we present a method for solving a composite optimization problem, where the smooth part of the objective is {\em quartic regular}. This means that its fourth-order derivative is uniformly positive definite and bounded on the domain of the composite term. We show that such a function admits global upper and lower bounds for their growth with respect to local approximations. Using the corresponding inequalities, we justify a linear rate of convergence of some variant of DQNM.

In Section \ref{sc-Conc}, we present some applications of our results. First of all, we show that DQNM can be applied to functions with positive-semidefinite fourth derivative. In terms of the number of calls of oracle, it has the rate of convergence $\tilde O(k^{-3})$. In the remaining part of the section we discuss applications of DQNM to the framework of proximal point methods. We show that it can be used for designing optimization methods with the rate of convergence $\tilde O(k^{-p})$ with $p \in \{ 3,4,5\}$. All these methods are of the order two.

\vspace{1ex}\noindent
{\bf Notation.}  
Let us fix a finite-dimensional real vector space $\E$. We measure distances in $\E$ by a general norm $\| \cdot \|$. Its dual space is denoted by $\E^{*}$. It is the space of all linear functions on $\E$, for which we define the norm in the standard way:
$$
\ba{rcl}
\| g \|^* & = & \max\limits_{x \in \E} \{ \; \la g, x \ra: \; \| x \| \leq 1 \; \}, \qquad g \in \E^*.
\ea
$$
Using this norm, we can define an induced norm for a self-adjoint linear operator $A: \E \to \E^*$ as follows:
$$
\ba{rcl}
\| A \| & = & \max\limits_{x \in \E} \{ | \la A x, x \ra | : \; \| x \| \leq 1 \}.
\ea
$$

For defining Euclidean norm, we introduce a positive-definite self-adjoint linear operator $B: \E \to \E^*$ (notation $B \succ 0$). Then for $x \in \E$ and $g \in \E^*$ we have
$$
\ba{rcl}
\| x \|_B & = & \la Bx, x \ra^{1/2}, \quad \| g \|_{B}^{*} \; = \; \la g, B^{-1} g \ra^{1/2}.
\ea
$$
If no ambiguity arise, we drop index $B$ in this notation. Note that for any self-adjoint linear operator $A: \E \to \E^*$ we have
$$
\ba{rcl}
\| A \| & = & \min \{ \lambda: \;  \lambda B \succeq \pm A \}.
\ea
$$
If in addition, $A$ is positive semidefinite and $L \geq \| A \|$, then
\beq\label{eq-MOrder}
\ba{rcl}
A B^{-1} A \; \preceq \; L A.
\ea
\eeq

For function $f(\cdot)$ with open domain $\dom f \subseteq \E$, we denote by $\nabla f(x) \in \E^*$ its gradient and by $\nabla^2 f(x): \E \to \E^*$ its Hessian at point $x \in \dom f$. For the mixed directional derivative of order $p \geq 1$ along directions $h_1, \dots, h_p$ in $\E$, we use notation $D^pf(x)[h_1, \dots, h_p]$. Note that under the standard continuity assumptions this is a symmetric $p$-linear form. If all directions are the same, we use notation $D^pf(x)[h]^p$.

We use the same notation for symmetric 4-linear forms. Let $f[x_1,x_2,x_3,x_4]$ be such a form with $x_i \in \E$, $i=1,\dots,4$. If all vectors are the same, we denote the value of this form by $f[x]^4$. This is a multivaried polynomial of degree four. Let us treat it as a function $\vf(x) = f[x]^4$. Then the vector $f[x]^3 \in \E^*$ is defined as follows:
$$
\ba{rcl}
f[x]^3 & = & {1 \over 4} \nabla \vf(x), \quad x \in \E.
\ea
$$
Note that $\la f[x]^3, x \ra = f[x]^4$. Similarly,
$$
\ba{rcl}
f[x]^2 & = & {1 \over 12} \nabla^2 \vf(x), \quad x \in \E.
\ea
$$
Then $f[x]^2[y]^2 = {1 \over 12} \la \nabla^2\vf(x) y, y \ra$ and $f[x]^2[x]^2 = f[x]^4$. Finally,
$$
\ba{rcl}
f[x] & = & {1 \over 24} D^3\vf(x), \quad x \in \E,
\ea
$$
and $f[x][y]^3 = {1 \over 24} D^3\vf(x) [y]^3$, $f[x][x]^3 = f[x]^4$. Thus, $f = {1 \over 24} D^4\vf$.

Finally, we need one result, which is usually proved under more restrictive assumptions. 
\BL\label{lm-GLip}
Let function $f(\cdot)$ be convex and twice continuously differentiable at its open domain. For two points $x, y \in \dom f$, denote $M_{x,y} = \max\limits_{0\leq \alpha \leq 1} \| \nabla^2 f(\alpha x + (1-\alpha)y) \|$, where the norm is Euclidean. Then
\beq\label{eq-GDiff}
\ba{rcl}
\| \nabla f(x) - \nabla f(y) \|_B^* & \leq & M_{x,y} \la \nabla f(x) - \nabla f(y), x - y \ra.
\ea
\eeq
\EL
\proof
Denote $G = \int\limits_0^1 \nabla^2 f(x + \tau (y - x)) d \tau$. Then $G(y-x) = \nabla f(y) - \nabla f(x)$ and $\| G \| \leq M_{x,y}$. Therefore, inequality (\ref{eq-GDiff}) follows from (\ref{eq-MOrder}).
\qed

\section{Positive semidefiite quartic forms}\label{sc-CQP}
\SetEQ

In the first part of this paper, we are interested in numerical methods for finding an approximate solution to the following minimization problem:
\beq\label{prob-Quart}
F^* \; = \; \min\limits_{x \in \dom \psi } \; \Big[ \; F(x) = f(x) + \psi(x) \; \Big],
\eeq
where $\psi(\cdot)$ is a {\em simple} closed convex function and $f(\cdot)$ is a convex {\em quartic polynomial}. 
Since, $f(\cdot)$ is a polynomial of degree four, its Taylor expansion with respect to any point in $\E$ is exact. Namely, for any $x, y \in \E$ we have
\beq\label{eq-Taylor}
\ba{rcl}
f(y) & = & f(x) + \la \nabla f(x), y-x] \ra + \frac{1}{2} \la \nabla^2 f(x)(y-x), y-x \ra  \\
\\
& & + \frac{1}{6} D^3f(x)[y-x]^3 + f_4[y-x]^4.
\ea
\eeq
where $f_4[\cdot,\cdot,\cdot, \cdot]$ is a symmetric four-linear form in $\E$. Thus, $f_4 = {1 \over 24} D^4f$.

Convexity of function $f(\cdot)$ has several important consequences.
\BT\label{th-Pos}
For any $x$ and $y$ in $\E$ we have 
\beq\label{eq-PosXY}
\ba{rcl}
0 &\leq & f_4[x]^2[y]^2.
\ea
\eeq
This property implies the following relations:
\beq\label{eq-Conv3}
\ba{rcl}
2 f_4[x,y] & \preceq & f_4[x]^2 + f_4[y]^2,
\ea
\eeq
\beq\label{eq-Cauchy1}
\ba{rcl}
(f_4[x]^2[y]^2)^2 & \leq & f_4[x]^4 \cdot f_4[y]^4,
\ea
\eeq
\beq\label{eq-Cauchy2}
\ba{rcl}
(f_4[x]^3[y])^2 & \leq & f_4[x]^4 \cdot f_4[x]^2[y]^2.
\ea
\eeq
\ET
\proof
Indeed, taking in (\ref{eq-Taylor}) $x = 0$, we have $\nabla^2 f(y) = \nabla^2 f(0) + D^3 f(0)[y] + \half D^4 f(0)[y]^2$. Therefore,
$$
\ba{rcl}
0 \; \leq \; \la \nabla^2 f(y) h, h \ra & = & \la \nabla^2 f(0) h, h \ra + D^3 f(0)[y][h]^2 + 12 f_4[y]^2[h]^2.
\ea
$$
Replacing now $y$ by $\tau y$ and taking the limit as $\tau \to \infty$, we get inequality (\ref{eq-PosXY}). It implies that $f_4[x-y]^2 \succeq 0$. And this is exactly the relation (\ref{eq-Conv3}). 

Further, in view of (\ref{eq-PosXY}), for any $\tau \in \R$ we have $f_4[x-\tau y]^2[x+\tau y]^2 \geq 0$. If we open the brackets, then we get
$$
\ba{rcl}
0 & \leq & f_4[x]^4 + \tau^4 f_4[y]^4 - 2 \tau^2 f_4[x]^2[y]^2.
\ea
$$
If $f_4[y]^4 = 0$, then $f_4[x]^2[y]^2 = 0$ and (\ref{eq-Cauchy1}) holds. If not, then by minimizing the right-hand side of the latter inequality in $\tau$, we get the relation (\ref{eq-Cauchy1}). Finally, minimizing the right-hand side of inequality 
$$
\ba{rcl}
0 & \leq & f_4[x]^2[x - \tau y ]^2
\ea
$$
in $\tau$, we get inequality (\ref{eq-Cauchy2}).
\qed

Define $d(x) = f_4[x]^4$. Note that
\beq\label{eq-Hom}
\ba{rcl}
d(tx) & = & |t|^4 d(x), \quad x \in \E, \; t \in \R.
\ea
\eeq
Let us look at some properties of this function. First of all, note that all derivatives of this function can be written in terms of the form $f_4$ as follows:
\beq\label{eq-DerD}
\ba{rclrcl}
\la \nabla d(x), h \ra & = & 4 f_4[x]^3[h], \quad \la \nabla^2 d(x) h, h \ra & = & 12 f_4[x]^2[h]^2,\\
\\
D^3d(x)[h]^3 & = & 24 f_4[x][h]^3, \quad D^4 d(x)[h]^4 & = & 24 f_4[h]^4.
\ea
\eeq
where $x$ and $h$ are arbitrary vectors in $\E$. Thus, by inequality (\ref{eq-PosXY}), we can see that function $d(\cdot)$ is convex.

Let us assume now that the form $f_4$ is positive definite:
\beq\label{ass-Pos}
\ba{rcl}
f_4[x]^4  & > & 0 \quad \forall x \in \E \setminus \{0\}.
\ea
\eeq
Then, we can define in $\E$ the following norm:
\beq\label{def-NormF}
\ba{rcl}
\| x \|_f & = & [d(x)]^{1/4}, \quad x \in \E.
\ea
\eeq

This norm has several interesting properties. First of all, it justifies the uniform convexity of function $d(\cdot)$ itself.
\BL\label{lm-DUni}
Function $d(\cdot)$ is uniformly convex of degree four with constant $\sigma_d = \frac{4}{3}$:
\beq\label{eq-DUni}
\ba{rcl}
d(y) & \geq & d(x) + \la d(x), y - x \ra + \frac{1}{4} \sigma_d \| y - x \|_f^4, \quad x, y \in \E.
\ea
\eeq
Moreover, we have the following lower quadratic bound for the growth of function $d(\cdot)$:
\beq\label{eq-DQuad}
\ba{rcl}
d(y) & \geq & d(x) + \la d(x), y - x \ra + \frac{1}{6} \la \nabla^2d(x)(y-x), y - x \ra, \quad x, y \in \E.
\ea
\eeq
\EL
\proof
Indeed, for any $x$ and $y$ from $\E$, we have
$$
\ba{c}
d(y) -  d(x) - \la d(x), y - x \ra \; = \; \half \la \nabla^2 d(x) (y-x), y- x \ra \\ \\
+ \frac{1}{6} D^3d(x)[y-x]^3 + \frac{1}{24} D^4d(x)[y-x]^4\\
\\
\refEQ{eq-DerD} \; 6 f_4[x]^2[y-x]^2 + 4 f_4[x][y-x]^3 + f_4[y-x]^4\\
\\
\refGE{eq-Cauchy2} \; 6 f_4[x]^2[y-x]^2 - 4 \left( f_4[x]^4 \cdot f_4[x]^2[y-x]^2 \right)^{1/2} + f_4[y-x]^4.
\ea
$$
Minimizing the right-hand side of this inequality in $f_4[x]^2[y-x]^2$, we get
$$
\ba{rcl}
d(y)- d(x) -  \la d(x), y - x \ra & \geq & (1 - \frac{2}{3}) f_4[y-x]^4 \; = \; \frac{1}{3} \| y - x \|^4_f.
\ea
$$
Finally, minimizing the same inequality in $f_4[x]^4$, we get
$$
\ba{rcl}
d(y)- d(x) -  \la d(x), y - x \ra & \geq & (6 - 4) f_4[x]^2[y-x]^2,
\ea
$$
and this is (\ref{eq-DQuad}).
\qed

Consider now the function $Q_f(x) = \half \| x \|^2_f$.
\BL\label{lm-Quad}
Function $Q_f(\cdot)$ is strongly convex and has Lipschitz continuous gradients 
with respect to the norm $\| \cdot \|_f$:
\beq\label{eq-QLip}
\ba{rcl}
\frac{1}{12 \| x \|_f^2} \la \nabla^2 d(x)h, h \ra & \leq & \la \nabla^2 Q_f(x)h,h \ra \; \leq \; 3 \| h \|^2_f, \quad x, h \in \E.
\ea
\eeq
\EL
\proof
Since $Q_f(x) = \half d^{1/2}(x)$, for any $x$ and $h$ from $\E$, we have
$$
\ba{rcl}
\la \nabla Q_f(x), h \ra & = & d^{-1/2}(x) f_4(x)[x]^3[h],\\
\\
\la \nabla^2 Q_f(x)h, h \ra & = & 3 d^{-1/2}(x) f_4(x)[x]^2[h]^2 - 2 d^{-3/2}(x) (f_4[x]^3[h])^2\\
\\
& \refGE{eq-Cauchy2} & d^{-1/2}(x) f_4(x)[x]^2[h]^2. 
\ea
$$
On the other hand,
$$
\ba{rcl}
\la \nabla^2 Q_f(x)h,, h \ra & \leq & 3 d^{-1/2}(x) f_4(x)[x]^2[h]^2 \; \refLE{eq-Cauchy1} \; 3 \| h \|_f^2. \QF
\ea
$$

\section{Minimizing convex polynomials of degree four}\label{sc-Min}
\SetEQ

In this section, we are going to check which properties of function $d(\cdot)$ are inherited by function $f(\cdot)$ and how we can use them for constructing minimization schemes. The next statement is a variant of Theorem 1 in \cite{ImpTensor}.
\BL\label{lm-SmallD3}
For any point $x$ and direction $h \in \E$, we have
\beq\label{eq-SmallD3}
\ba{rcl}
D^3f(x)[h] & \preceq & {1 \over \tau} \nabla^2 f(x) + 12 \tau f_4[h]^2
\ea
\eeq
\EL
\proof
Since $f(\cdot)$ is a convex quartic polynomial, for any $\tau > 0$ we have
$$
\ba{rcl}
0 & \preceq & \nabla^2f(x - \tau h) \; = \; 
\nabla^2f(x) - \tau D^3f(x)[h] + {\tau^2 \over 2} D^4f(x)[h]^2.
\ea
$$
Dividing this inequality by $\tau$, we get
$$
\ba{rcl}
D^3f(x)[h] & \preceq & {1 \over \tau} \nabla^2f(x) +{\tau \over 2} D^4 f(x)[h]^2,
\ea
$$
and this is the relation (\ref{eq-SmallD3}).
\qed

\BC\label{cor-ABounds}
For all $x, y \in \E$ and $\tau > 0$, we have
\beq\label{eq-UBound}
\ba{rl}
& f(y) - f(x) - \la f(x), y - x \ra - \half \la \nabla^2 f(x) (y-x), y - x \ra \\
\\
\leq & {1 \over 6 \tau} \la \nabla^2 f(x)(y-x), y - dx \ra + (1+2\tau) \| y - x \|^4_f,
\ea
\eeq
\beq\label{eq-LBound}
\ba{rl}
& f(y) - f(x) - \la f(x), y - x \ra - \half \la \nabla^2 f(x) (y-x), y - x \ra \\
\\
\geq & - {1 \over 6 \tau} \la \nabla^2 f(x)(y-x), y - dx \ra + (1-2\tau) \| y - x \|^4_f.
\ea
\eeq
\EC
\proof
In view of relation (\ref{eq-SmallD3}), we have
$| D^3f(x)[h]^3 |\leq {1 \over \tau} \la \nabla^2f(x)h,h \ra + 12 \tau \| h \|^4_f$.
Therefore, 
$$
\ba{rl}
& f(y) - f(x) - \la f(x), y - x \ra - \half \la \nabla^2 f(x) (y-x), y - x \ra \\
\\
\refEQ{eq-Taylor} & {1 \over 6} D^3f(x)[y-x]^3 +  \| y - x \|^4_f \; \leq \; {1 \over 6 \tau} \la \nabla^2f(x)h,h \ra + (1 + 2\tau ) \| h \|^4_f. 
\ea
$$
The second inequality can be justified in the same way.
\qed

Using inequalities (\ref{eq-UBound}) and (\ref{eq-LBound}), we can write down the lower and upper approximations for the smooth part of the objective function:
\beq \label{eq-Apps}
\ba{rcl}
\bar f_{x,\tau}(y) & = & f(x) + \la f(x), y - x \ra + {3\tau - 1 \over 6 \tau} \la \nabla^2 f(x) (y-x), y - x \ra \\
\\
& & + (1- 2\tau) \| y - x \|^4_f \; \leq \; f(y),\\
\\
\hat f_{x,\tau}(y) & = & f(x) + \la f(x), y - x \ra + {1 + 3\tau \over 6 \tau} \la \nabla^2 f(x) (y-x), y - x \ra \\
\\
& & + (1+2\tau ) \| y - x \|^4_f \; \geq \; f(y), \quad x, y \in \E.
\ea
\eeq
Note that we are interested in $\tau \in \left[{1 \over 3}, \half \right]$.

Let us assume that the norm $\| \cdot \|_f$ is simple enough for minimizing the following auxiliary function,
$$
\hat f_{x,\tau}(y) + \psi(y), \quad y \in \dom \psi,
$$
in a closed form. Thus, we can use it for justifying methods for solving the problem (\ref{prob-Quart}).

Let us define the following parametric function:
\beq\label{def-FPar}
\ba{rcl}
\bar f_{x,\tau}(\alpha,y) & = & f(x) + \la f(x), y - x \ra + {3\tau - 1 \over 6 \tau \alpha } \la \nabla^2 f(x) (y-x), y - x \ra\\
\\
& & + {1- 2\tau \over \alpha^3} \| y - x \|^4_f, \quad y \in \E, \; \alpha > 0.
\ea
\eeq
Note that $\bar f_{x,\tau}(1,y) = \bar f_{x,\tau}(y) \leq f(y)$ for all $y \in \E$. On the other hand, function $\bar f_{x,\tau}(\cdot,\cdot)$ is jointly convex in its arguments. Therefore, the univariate function
$$
\ba{rcl}
F^*_{x,\tau}(\alpha) & = & \min\limits_{y \in \dom \psi} \Big[ \; \bar f_{x,\tau}(\alpha,y) + \psi(y) \; \Big], \quad \alpha > 0,
\ea
$$
is convex. Let us prove the following result.
\BT\label{th-Alpha}
Let $\tau \in \left({1 \over 3}, {1 \over 2} \right)$ and the parameter $\alpha > 0$ satisfies the following condition:
\beq\label{eq-Alpha}
\ba{rcl}
\alpha & \leq & \min \Big\{ {3 \tau - 1 \over 3 \tau + 1}, \left[ 1 - 2 \tau \over 1 + 2 \tau \right]^{1/3} \Big\}.
\ea
\eeq
Then
\beq\label{eq-DecF}
\ba{rcl}
\min\limits_{y \in \dom \psi} \Big[ \; \hat f_{x,\tau}(y) + \psi(y) \; \Big] & \leq & \alpha F^* + (1-\alpha) F(x).
\ea
\eeq
\ET
\proof
Note that $F^*_{x,\tau}(1) \leq F^*$. On the other hand, by continuity, we can define $F^*_{x,\tau}(0) = F(x)$. At the same time, condition (\ref{eq-Alpha}) ensures the following relations:
$$
\ba{rcl}
{3 \tau - 1 \over 6 \tau \alpha} & \geq & {3 \tau + 1 \over 6 \tau }, \quad {1 - 2 \tau \over \alpha^3} \; \geq \; 1 + 2 \tau.
\ea
$$
This means that $\bar f_{x,\tau}(\alpha,y) \geq \; \hat f_{x,\tau}(y)$ for all $y \in \E$. Hence,
$$
\ba{rcl}
\min\limits_{y \in \dom \psi} \Big[ \;  \hat f_{x,\tau}(y) + \psi(y) \; \Big] & \leq & F^*_{x,\tau}(\alpha) \; \leq \; 
\alpha F^*_{x,\tau}(1) + (1 - \alpha) F^*_{x,\tau}(0) \\
\\
& \leq & \alpha F^* + (1-\alpha) F(x). \QR
\ea
$$

Inequality (\ref{eq-DecF}) is important for justification of the following optimization scheme. Denote by $\tau_*$ the unique positive root of the following equation:
\beq\label{eq-Tau}
\ba{rcl}
{3 \tau - 1 \over 3 \tau + 1} & = & \left[ 1 - 2 \tau \over 1 + 2 \tau \right]^{1/3} .
\ea
\eeq 
\beq\label{met-QN}
\ba{|c|}
\hline \\
\mbox{\bf Damped Quartic Newton Method}\\
\\
\hline \\
\ba{l}
\mbox{{\bf Initialization.} Choose $x_0 \in \dom \psi$.}\\
\\
\mbox{{\bf $k$th iteration ($k \geq 0$).} Iterate $x_{k+1} = \arg\min\limits_{y \in \dom \psi} \Big[ \; \hat f_{x_k, \tau_*}(y) + \psi(y) \; \Big]$.}\\
\\
\hline 
\ea
\ea
\eeq
\BT\label{th-Quartic}
For all $k \geq 0$, we have
\beq\label{eq-Quartic}
\ba{rcl}
F(x_k) - F^* & \leq & (1 - \alpha_*)^k(F(x_0) - F^*),
\ea
\eeq
where $\alpha^* = {3 \tau_*-1 \over 3 \tau_*+1} > 0.193$ with $\tau_* = {1 \over 6}\sqrt{3 + \sqrt{33}}$.
\ET
\proof
Note that for any $x,y \in \E$ we have $\bar f_{x,\tau_*}(\alpha_*,y) = \hat f_{x,\tau_*}(y)$. Therefore,
$$
\ba{rcl}
F(x_{k+1}) & \refLE{eq-Apps} & \hat f_{x_k,\tau_*}(x_{k+1}) + \psi(x_{k+1}) \; \refLE{eq-DecF}\; \alpha_* F^* + (1- \alpha_*) F(x_k).
\ea
$$
It remains to note that $\tau_*$ is the root of the equation
$$
\ba{rcl}
(3 \tau-1)^3 (1+2\tau) & \refEQ{eq-Tau} & (3 \tau +1)^3 (1-2\tau).
\ea
$$
After simplification, it is reduced to $54 \tau^4 - 9 \tau^2 = 1$, which gives us $\tau_* = {1 \over 6}\sqrt{3 + \sqrt{33}}$.
\qed

Method (\ref{met-QN}) is implementable if the structure of the norm $\| \cdot \|_f$ is quite simple. Let us end up this section with an important source of such problems.
\BE\label{ex-T3}
For a convex function, consider its Taylor polynomials of degree three, augmented by a fourth power of Euclidean norm:
$$
\ba{rcl}
f_{\bar x, H}(y) & = & f(\bar x) + \la \nabla f(\bar x), y - \bar x \ra + \half \la \nabla^2f(\bar x)(y- \bar x), y - \bar x\ra\\
\\
& & + {1 \over 6} D^3f(\bar x)[y-\bar x]^3 + {H \over 24} \| y - \bar x \|^4,
\ea
$$
where $\| \cdot \|$ is a Euclidean norm.
The possibility to minimize this polynomial is essential for the third-order schemes. It can be proved \cite{ImpTensor} that for $H \geq 3 L_3$, where $L_3$ is the Lipschitz constant for third derivative of function $f(\cdot)$, the quartic polynomial $f_{\bar x, H}(\cdot)$ is convex. Note that it fits our framework since $D^4 f_{\bar x, H}(x)[h]^4 = H \| h \|^4$. Thus, in this case, the auxiliary problem in method (\ref{met-QN}) can be solved by using the standard tools of Linear Algebra.
\EE

\section{Simplified approximations}\label{sc-Simple}
\SetEQ

Let us introduce in the space $\E$ some norm $\| \cdot \|$, which is simpler than the norm $\| \cdot \|_f$.
Since the space $\E$ is finite-dimensional, assumption (\ref{ass-Pos}) implies existence of two constants $0 < \mu \leq L$ such that
\beq\label{eq-PhiBound}
\ba{rcl}
\mu \| x \|^4 & \leq & f_4[x]^4 \; \leq \; L \| x \|^4, \quad x \in \E.
\ea
\eeq
In the particular case when $\| \cdot \| = \| \cdot \|_f$, we have $\mu = L = 1$. Denote by $q = {\mu \over L}$ the condition number of the norm $\| \cdot\|$ with respect to the optimal norm $\| \cdot \|_f$.

Now, in view of (\ref{eq-Apps}), we can define the following simplified lower and upper bounds for the smooth part of the objective function:
\beq \label{eq-Apps1}
\ba{rcl}
\bar \phi_{x,\tau}(y) & \Def & f(x) + \la f(x), y - x \ra + {3\tau - 1 \over 6 \tau} \la \nabla^2 f(x) (y-x), y - x \ra \\
\\
& & + (1- 2\tau) \mu \| y - x \|^4 \; \leq \; f(y),\\
\\
\hat \phi_{x,\tau}(y) & \Def & f(x) + \la f(x), y - x \ra + {1 + 3\tau \over 6 \tau} \la \nabla^2 f(x) (y-x), y - x \ra \\
\\
& & + (1+2\tau ) L \| y - x \|^4 \; \geq \; f(y), \quad x, y \in \E.
\ea
\eeq
As in (\ref{def-FPar}), define the following parametric functions
\beq\label{def-PhiPar}
\ba{rcl}
\bar \phi_{x,\tau}(\alpha, y) & \Def & f(x) + \la f(x), y - x \ra + {3\tau - 1 \over 6 \tau \alpha} \cdot \la \nabla^2 f(x) (y-x), y - x \ra \\
\\
& & + (1- 2\tau) \cdot {\mu \over \alpha^3} \| y - x \|^4,\\
\\
\Phi^*_{x,\tau}(\alpha) & = & \min\limits_{y \in \dom \psi} \Big[ \; \bar \phi_{x,\tau}(\alpha, y) + \psi(y) \; \Big], \quad \alpha > 0.
\ea
\eeq

Consider the following optimization method.
\beq\label{met-QN2}
\ba{|c|}
\hline \\
\mbox{\bf Relaxed Quartic Newton Method}\\
\\
\hline \\
\ba{l}
\mbox{{\bf Initialization.} Choose $x_0 \in \dom \psi$ and parameter $\tau \in \left({1 \over 3}, \half \right)$.}\\
\\
\mbox{{\bf $k$th iteration ($k \geq 0$).} Iterate $x_{k+1} = \arg\min\limits_{y \in \dom \psi} \Big[ \; \hat \phi_{x_k, \tau}(y) + \psi(y) \; \Big]$.}\\
\\
\hline 
\ea
\ea
\eeq
As compared with metod (\ref{met-QN}), we use here a conservative approximation of the smooth part of the objective function: $\hat \phi_{x,\tau}(y) \refGE{eq-PhiBound} \hat f_{x,\tau}(y)$ for all $x, y \in \E$.
\BT\label{th-Quartic2}
Let the parameter $\tau$ of method (\ref{met-QN2}) satisfy the following condition:
\beq\label{eq-Tau2}
\ba{rcl}
\alpha & \leq & \min \Big\{ {3 \tau - 1 \over 3 \tau + 1} , \left[ q \cdot { 1 - 2 \tau \over 1 + 2 \tau} \right]^{1/3} \Big\}
\ea
\eeq 
for some $\alpha > 0$.
Then, for all $k \geq 0$, we have
\beq\label{eq-Quartic2}
\ba{rcl}
F(x_k) - F^* & \leq & (1 - \alpha)^k(F(x_0) - F^*).
\ea
\eeq
\ET
\proof
Indeed, in view of condition (\ref{eq-Tau2}), we have
$$
\ba{rcl}
{3 \tau + 1 \over 6 \tau} & \leq & {3 \tau - 1 \over 6 \tau \alpha}, \quad
(1 + 2 \tau) L \; \leq \; (1 - 2\tau ) \cdot { \mu \over \alpha^3}.
\ea
$$
Hence, 
$$
\ba{rcl}
F(x_{k+1}) & \refLE{eq-Apps1} & \hat \phi_{x_k,\tau}(x_{k+1}) + \psi(x_{k+1}) \; \refEQ{met-QN2} \; \min\limits_{y \in \dom \psi} \Big[ \;
\hat \phi_{x_k,\tau}(y) + \psi(y)\; \Big]\\
\\
& \refLE{eq-Tau2} & \min\limits_{y \in \dom \psi} \Big[ \;
\bar \phi_{x_k,\tau}(\alpha, y) + \psi(y)\; \Big] \; \refEQ{def-PhiPar} \; \Phi^*_{x_k,\tau}(\alpha).
\ea
$$
On the other hand, since function $\Phi^*_{x_k,\tau}(\cdot)$ is convex, we have
$$
\ba{rcl}
\Phi^*_{x_k,\tau}(\alpha) & \leq & \alpha \Phi^*_{x_k,\tau}(1) + (1 - \alpha) \Phi^*_{x_k,\tau}(0) \; \refLE{eq-Apps1} \; \alpha F^* + (1 - \alpha) F(x_k).
\ea
$$
Thus, we obtain the bound (\ref{eq-Quartic2}).
\qed

Since we are interested in the maximal possible value of $\alpha$, let us try to understand which value of $\tau \in \left({1 \over 3}, {1 \over 2} \right)$ makes the right-hand side of inequality (\ref{eq-Tau2}) big. Note that 
$$
\ba{rcl}
{3 \tau - 1 \over 3 \tau + 1} & \geq & {2 \over 5} (3 \tau -1), \quad \tau \in \left({1 \over 3}, {1 \over 2} \right).
\ea
$$
Let us look now at the function $\xi(\tau) = \omega^{1/3}(\tau)$ with $\omega(\tau) = {1 - 2 \tau \over 1 + 2 \tau}$. Since
$$
\ba{rcl}
\xi'(\tau) & = & {1 \over 3} \omega^{-2/3}(\tau) \omega'(\tau) \quad \xi''(\tau) = -{2 \over 9} \omega^{-5/3}(\tau) (\omega'(\tau))^2 + {1 \over 3} \omega^{-2/3} \omega''(\tau),\\
\\
\omega'(\tau) & = & - {4 \over (1+2\tau)^2}, \quad \omega''(\tau) \; = \; {16 \over (1+2 \tau)^3},
\ea
$$
we conclude that $\xi''(\tau) = {16 \omega^{-5/3}(\tau) \over 9 (1+2 \tau)^4} (- 2 + 3(1-2 \tau) ) \leq 0$ for $\tau \in \left({1 \over 3}, {1 \over 2} \right)$. Hence, for this interval, we have
$$
\ba{rcl}
\xi(\tau) & \geq & \xi({1 \over 3}) \cdot 3(1-2 \tau) \; = \; 3(1-2 \tau) \cdot 5^{-1/3}  .
\ea
$$
Thus, the condition (\ref{eq-Tau2}) is satisfied for 
\beq\label{eq-Alpha2}
\ba{rcl}
\alpha & = & \min\Big\{ {2 \over 5} (3 \tau -1), 3(1-2 \tau) \cdot \left( q \over 5 \right)^{1/3} \Big\}.
\ea
\eeq
Therefore, the best choice of $\tau$ is as follows:
\beq\label{eq-TOpt2}
\ba{rcl}
\tau_{\#} & = & \half - {1 \over 6 (1+5\kappa)}, \quad \kappa \; = \; \left( q \over 5 \right)^{1/3}.
\ea
\eeq
In this case, $\alpha_{\#} = {2 \over 5} (3 \tau_{\#} -1) \; = \; {\kappa \over 1 + 5 \kappa}$. Thus, we have proved the following theorem.
\BT\label{th-RQN}
Let the sequence of points $\{ x_k \}_{k \geq 0}$ be generated by the method (\ref{met-QN2}) with parameter $\tau = \tau_{\#}$. Then, for all $k \geq 0$, we have
\beq\label{eq-RateRQN}
\ba{rcl}
F(x_k) - F^* & \leq & (1 - \alpha_{\#})^k(F(x_0) - F^*).
\ea
\eeq
\ET

In order to generate an $\epsilon$-solution of problem (\ref{prob-Quart})   in function value, method~(\ref{met-QN2}) needs at most
\beq\label{eq-CompRQN}
\ba{c}
\left(5 + \left({5 \over q}\right)^{1/3} \right) \ln {F(x_0) - F^* \over \epsilon}
\ea
\eeq
iterations. Note that for $q = 1$, we have $\alpha_{\#} = {1 \over 5 + 5^{1/3}} > 0.149$, which is quite close to $\alpha_*$. Thus, the general rule (\ref{eq-TOpt2}) is not much worse than the optimal specific rule (\ref{eq-Tau}). 

For practical applications, it is reasonable to use the Euclidean norm $\| \cdot \|$. In this case, the complexity of one iteration of method (\ref{met-QN2}) is the same as that of the usual Cubic Newton Method \cite{CNM}

\section{Quartic regularity for convex functions}\label{sc-QReg}
\SetEQ

Let us check to what extent the results of the previous sections can be propagated onto the class of general convex functions. Consider the convex unconstrained minimization problem
\beq\label{prob-QRMin}
\ba{rcl}
F^* & = & \min\limits_{x \in \dom \psi} \Big[ \; F(x) \; = \; f(x) + \psi(x) \; \Big],
\ea
\eeq
where $\psi(\cdot)$ is a simple closed convex function, and the smooth part $f(\cdot)$ of the objective satisfies an additional assumption of {\em quartic regularity}.
\BD\label{DEF-QReg}
Function $f(\cdot)$ with convex and open domain $\dom f \subseteq \E$ is called {\em quartic regular} (Q-regular) on convex set $S$, if it is four times continuously differentiable on $\dom f$, and there exist two parameters $L \geq \mu \geq 0$ such that
\beq\label{def-QRegL}
\ba{rcl}
D^4f(x)[h]^4 \; \geq \; \mu \| h \|^4, \quad \forall h \in \E, \; x \in S, 
\ea
\eeq
\beq\label{def-QRegU}
\ba{rcl}
D^4f(x)[h]^4 \; \leq \; L \| h \|^4, \quad \forall h \in \E, \; x \in S_{1/2}, 
\ea
\eeq
where $S_{1/2} = \{ z = x - \half (y-x): \; x, y \in S \} \subseteq \dom f$.
\ED
This definition can be applied to general norms.\footnote{$^)$ As compared to (\ref{eq-PhiBound}), it is more convenient now to work directly with the fourth derivative.}$^)$

Note that quartic regularity itself does not imply convexity. Thus, for this section, our assumption is as follows.
\BAS\label{ass-QR}
Smooth part of the objective function in problem (\ref{prob-QRMin}) is quartic regular on the set $S \supseteq \dom \psi$.
\EAS

Let us mention some trivial properties of Q-regular functions.
\BI
\II
If functions $f_i(\cdot)$ are Q-regular on sets $S_i$ with parameters $(\mu_i,L_i)$, $i=1,2$, then their weighted sum $f(x) = \alpha f_1(x) + \beta f_2(x)$, $\alpha, \beta \geq 0$, is Q-regular on $S_1 \bigcap S_2$ with parameters 
$$
(\alpha \mu_1 + \beta \mu_2, \alpha L_1 + \beta L_2).
$$
\II
If function $\phi(\cdot)$ is Q-regular on $S_1$ with parameters $(\mu,L)$ and $A: \E \to \E_1$ is a non-degenerate linear operator, then function $f(x) = \phi(Ax)$ is Q-regular on $S = \{x: \; A x \in S_1 \}$ with parameters 
$(\mu \, \sigma^4_{\min}(A), L \, \sigma^4_{\max}(A))$.
\II
In particular, if function $f(\cdot)$ is Q-regular on $S$ with parameters $(\mu,L)$, then function $\vf_{\tau}(x) = f(\tau x)$, with scaling factor $\tau\neq 0$, is Q-regular on $\tau^{-1}Q$ with parameters $(\mu \tau^4, L \tau^4)$.
\II
Function $d(x) = {1 \over 24} \| x \|^4$, defined by a Euclidean norm, is Q-regular on $\E$ with parameters 
$\mu = L = 1$.
\EI

The abilities of optimization methods in solving the problem (\ref{prob-QRMin}) are supported by the following inequalities.
\BT\label{th-QR}
Let function $f(\cdot)$ be Q-regular on $S$ with parameters $\mu$ and $L$. Then for all $x, y \in S$ and $\gamma \in \left[ {1 \over 3}, {1 \over 2} \right]$, we have
\beq\label{eq-QReg}
\ba{rcl}
f(y) & \geq & f(x) + \la \nabla f(x), y - x \ra + {3 \gamma - 1 \over 6 \gamma} \la \nabla^2 f(x)(y-x), y - x \ra\\
\\
& &  +{\mu \over 24} \| y - x \|^4 \Big[ (1-2 \gamma) - 
{16 \over 125}  \left({L \over \mu}-1\right)\left({3 \gamma - 1 \over \gamma} \right)^3 \Big].
\ea
\eeq
Moreover, for all $\gamma > 0$, we have
\beq\label{eq-QRegUp}
\ba{rcl}
f(y) & \leq & f(x) + \la \nabla f(x), y - x \ra + {3 \gamma + 1 \over 6 \gamma} \la \nabla^2 f(x)(y-x), y - x \ra\\
\\
& &  + {2 \gamma + 1 \over 24} L \| y - x \|^4.
\ea
\eeq
\ET
\proof
By Taylor formula, we have
$$
\ba{rcl}
\Delta(x,y) & \Def & f(y) - f(x) - \la \nabla f(x), y - x \ra - \half \la \nabla^2 f(x)(y-x), y - x \ra\\
& = & {1 \over 6} D^3f(x)[y-x]^3 + {1 \over 3!} \int\limits_0^1 (1-\tau)^3 D^4f(x+\tau(y-x)) [y-x]^4 d \tau.
\ea
$$
Since $f(\cdot)$ is convex,
$$
\ba{rcl}
0 \; \preceq \; \nabla^2 f(x + \gamma (y-x)) & = & \nabla^2 f(x) + \gamma D^3f(x)(y - x) \\
\\
& &  + \int\limits_0^\gamma (\gamma - \tau) D^4 f(x+\tau (y-x))[y-x]^2 d \tau.
\ea
$$
Hence,
\beq\label{eq-TLow}
\ba{rcl}
\Delta(x,y) & \geq & - {1 \over 6 \gamma} \la \nabla^2 f(x)(y-x), y - x \ra \\
\\
& & + \int\limits_0^{\gamma} \left[ {(1 - \tau)^3 \over 3!} - {1 \over 6}\left(1 - {\tau \over \gamma} \right) \right] D^4 f(x+\tau (y-x))[y-x]^4 d \tau\\
\\
& & + {1 \over 3!} \int_{\gamma}^1 (1-\tau)^3 D^4f(x+\tau(y-x)) [y-x]^4 d \tau.
\ea
\eeq	
Denote $\theta(\tau) = (1-\tau)^3 + {\tau \over \gamma} - 1 = \tau \left( {1 \over \gamma} - 3 \right) + 3 \tau^2 - \tau^3$. This cubic polynomial has three real roots, one at $\tau = 0$ and two others at the points
$$
\ba{rcl}
\tau_{0,1}(\gamma) & = & {3 \over 2} \pm \sqrt{{1 \over \gamma}- {3 \over 4}}.
\ea
$$
Note that by  conditions of the theorem, we have
$$
\ba{rcl}
0 & \leq & \tau_0(\gamma) \; \leq \; \gamma.
\ea
$$
Therefore, at the interval $\tau \in [0,1]$, we have the following relations:
$$
\ba{rcl}
\theta(\tau) & \leq & 0, \quad \tau \in \Big[0,\tau_{0}(\gamma) \Big], \quad \theta(\tau) \; \geq \; 0, \quad  \tau_0(\tau) \leq \tau \leq 1 \quad ( \leq \tau_1(\gamma)).
\ea
$$
Denoting the first integral in (\ref{eq-TLow}) by $I_1$, the second one by $I_2$, and $r = \| y - x \|$, we have
$$
\ba{c}
6 I_1 \; = \; \int\limits_0^{\tau_0(\gamma)} \theta(\tau) D^4 f(x+\tau (y-x))[y-x]^4 d \tau + \int\limits_{\tau_0(\gamma)}^{\gamma} \theta(\tau) D^4 f(x+\tau (y-x))[y-x]^4 d \tau\\
\\
\refGE{def-QRegU} \; L r^4 \int\limits_0^{\tau_0(\gamma)} \theta(\tau) d \tau + \mu r^4  \int\limits_{\tau_0(\gamma)}^{\gamma} \theta(\tau)  d \tau\ \\
\\
= \; (L-\mu) r^4 \int\limits_0^{\tau_0(\gamma)} \theta(\tau) d \tau + \mu r^4  \int\limits_{0}^{\gamma} \theta(\tau)  d \tau .
\ea
$$
On the other hand, $6 I_2 \refGE{def-QRegL}  \mu r^4 \int_{\gamma}^1 (1-\tau)^3 d \tau$. Note that
$$
\ba{rcl}
\int\limits_{0}^{\gamma} \theta(\tau)  d \tau + \int_{\gamma}^1 (1-\tau)^3 d \tau & = & \int_{0}^1 (1-\tau)^3 d \tau + \int\limits_0^{\gamma} ({\tau \over \gamma} - 1) d \tau \; = {1 - 2 \gamma \over 4},
\ea
$$
and we conclude that
$$
\ba{rcl}
6(I_1+I_2) & \geq & (L-\mu) r^4 \int\limits_0^{\tau_0(\gamma)} \theta(\tau) d \tau + \mu r^4 {1 - 2 \gamma \over 4}.
\ea
$$
It remains to estimate the integral in the latter inequality from below.

Denoting for short $\tau_0 = \tau_0(\gamma)$, we have
$$
\ba{rcl}
I_0 & \Def & \int\limits_0^{\tau_0} \theta(\tau) d \tau \; = \; \int\limits_0^{\tau_0} \left[ (1-\tau)^3 + {\tau \over \gamma} - 1 \right] d \tau\\
\\
& = & - {1 \over 4}(1-\tau)^4 + {\tau^2 \over 2 \gamma} - \tau \Big|_0^{\tau_0} \; = \;  - {1 \over 4}(1-\tau_0)^4 + {\tau_0^2 \over 2 \gamma} - \tau_0 + {1 \over 4}
\ea
$$
Since $\theta(\tau_0) = 0$, we conclude that
$$
\ba{rcl}
I_0 & = & - {1 \over 4}(1-\tau_0)\left(1 - {\tau_0 \over \gamma} \right) + {\tau_0^2 \over 2 \gamma} - \tau_0 + {1 \over 4} \; = \; - {3 \gamma - 1 \over 4 \gamma} \tau_0 + {\tau_0^2 \over 4 \gamma}.
\ea
$$
Denote $\xi = {3 \gamma - 1 \over \gamma}$ and $\lambda = {4 \over 9} \xi$. Then $0 \leq \xi \leq 1$, $0 \leq \lambda \leq {4 \over 9}$, and
$$
\ba{c}
\tau_0 = {3 \over 2} - \sqrt{{1 \over \gamma} - {3 \over  4}} = {3 \over 2} - \sqrt{{9 \over 4} - \xi} = {3 \over 2} (1 - \sqrt{1 - \lambda}) \; \leq \; {3 \over 2} \cdot {3 \over 5} \lambda \; = \; {2 \over 5} \xi.
\ea
$$
Note that in our new notation we have
\beq\label{eq-Xi}
\ba{rcl}
3 \tau_0 - \tau_0^2 & = & 3 - {1 \over \gamma} \; = \; \xi.
\ea
\eeq
Therefore,
$$
\ba{rcl}
4 I_0 & = & - \xi \tau_0 + {1 \over \gamma} \tau_0^2 \; \refEQ{eq-Xi} \; -(3 \tau_0 - \tau_0^2)\tau_0 + {1 \over \gamma} \tau_0^2 \; = \; \tau_0^3 - \xi \tau_0^2\\
\\
& \refEQ{eq-Xi} & \tau_0^3 - (3 \tau_0 - \tau_0^2) \tau_0^2 \; = \; - 2 \tau_0^3 + \tau_0^4 \; \geq \; - 2 \tau_0^3 \; \geq \; -{16 \over 125} \xi^3.
\ea
$$

Thus, we get the following lower bound:
$$
\ba{rcl}
f(y) - f(x) - \la \nabla f(x), y - x \ra & \geq & {\xi \over 6} \la \nabla^2 f(x)(y-x), y - x \ra \\
\\
& & + {\mu r^4 \over 6} \Big[ {1 - 2 \gamma \over 4} -  {4 \over 125}  \left({L \over \mu}-1\right)\xi^3 \Big].
\ea
$$

In order to prove the upper bound, note that 
$$
\ba{rcl}
0 \; \preceq \; \nabla^2 f(x - \gamma (y-x)) & = & \nabla^2 f(x) - \gamma D^3f(x)(y - x) \\
\\
& &  + \int\limits_0^\gamma (\gamma - \tau) D^4 f(x-\tau (y-x))[y-x]^2 d \tau.
\ea
$$
Hence,
\beq\label{eq-TUp}
\ba{rcl}
\Delta(x,y) & \leq & {1 \over 6 \gamma} \la \nabla^2 f(x)(y-x), y - x \ra \\
\\
& & + {1 \over 6 \gamma} \int\limits_0^{\gamma}(\gamma - \tau) D^4 f(x-\tau (y-x))[y-x]^4 d \tau\\
\\
& & + {1 \over 3!} \int_{0}^1 (1-\tau)^3 D^4f(x+\tau(y-x)) [y-x]^4 d \tau\\
\\
& \refLE{def-QRegU} & {1 \over 6 \gamma} \la \nabla^2 f(x)(y-x), y - x \ra + {L \over 6} \| y - x \|^4 \left( { \gamma \over 2} + {1 \over 4} \right). \QF
\ea
\eeq	

Now we can estimate performance of method (\ref{met-QN2}) as applied to problem (\ref{prob-QRMin}). Define
\beq\label{def-Psi}
\ba{rcl}
\hat \xi_{x,\gamma}(y) & \Def & f(x) + \la f(x), y - x \ra + {1 + 3\gamma \over 6 \gamma} \la \nabla^2 f(x) (y-x), y - x \ra \\
\\
& & + {1+2\gamma \over 24} L \| y - x \|^4 \; \refGE{eq-QRegUp} \; f(y), \quad x, y \in S,
\ea
\eeq
where $f(\cdot)$ is a Q-regular function. 
Consider the following optimization method.
\beq\label{met-QN3}
\ba{|c|}
\hline \\
\mbox{\bf Quartic Regularization of Newton Method}\\
\mbox{\bf for Q-Regular Functions}\\
\\
\hline \\
\ba{l}
\mbox{{\bf Initialization.} Choose $x_0 \in \dom \psi$ and parameter $\gamma > 0$.}\\
\\
\mbox{{\bf $k$th iteration ($k \geq 0$).} Iterate $x_{k+1} = \arg\min\limits_{y \in \dom \psi} \Big[ \; \hat \xi_{x_k, \gamma}(y) + \psi(y) \; \Big]$.}\\
\\
\hline 
\ea
\ea
\eeq
\BT\label{th-Quartic3}
Let function $f(\cdot)$ be Q-regular on the sets $S \supseteq \dom \psi$ with parameters $0 < \mu \leq L$. And let the parameter $\gamma \in \left({1 \over 3}, \half \right)$ of method (\ref{met-QN3}) satisfy the following condition:
\beq\label{eq-Gamma3}
\ba{rcl}
\kappa(\gamma) & \Def & (1-2 \gamma) - {16 \over 125}  \left({L \over \mu}-1\right)\left({3 \gamma - 1 \over \gamma} \right)^3 \; \geq \; 0,\\
\\
\alpha & \leq & \min \Big\{ {3 \gamma - 1 \over 3 \gamma + 1} , 
\left[{\mu \over L(1+2\gamma)} \kappa(\gamma) \right]^{1/3} \Big\}
\ea
\eeq 
for some $\alpha > 0$.
Then, for all $k \geq 0$, we have
\beq\label{eq-Quartic3}
\ba{rcl}
F(x_k) - F^* & \leq & (1 - \alpha)^k(F(x_0) - F^*).
\ea
\eeq
\ET
\proof
As in the proof of Theorem \ref{th-Quartic2}, define the following parametric functions:
$$
\ba{rcl}
\hat \xi_{x,\gamma}(\alpha,y) & = & f(x) + \la \nabla f(x), y - x \ra + {3 \gamma - 1 \over 6 \alpha \gamma} \la \nabla^2 f(x)(y-x), y - x \ra \\
\\
& & +{\mu \over 24 \alpha^3} \| y - x \|^4 \kappa(\gamma),\\
\\
\Xi^*_{x,\gamma}(\alpha) & = & \min\limits_{y \in \dom \psi} \Big[ \; \hat \xi_{x,\gamma}(\alpha,y) + \psi(y) \; \Big], \quad \alpha > 0.
\ea
$$
Note that $\Xi^*_{x,\gamma}(0) = F(x)$ and $ \hat \xi_{x,\gamma}(1,y) \refLE{eq-QReg} f(y)$. Thus, $\Xi^*_{x,\gamma}(1) \leq F^*$. On the other hand, in view of conditions (\ref{eq-Gamma3}), we have
$$
\ba{rcl}
{3 \gamma - 1 \over 6 \alpha \gamma} & \geq & {3 \gamma +1 \over 6 \gamma}, \quad {\mu \over 24 \alpha^3} \kappa(\gamma) \; \geq \; {1 + 2 \gamma \over 24}L.
\ea
$$
Therefore, as in the proof of Theorem \ref{th-Quartic2}, we conclude that
$$
\ba{rcl}
F(x_{k+1}) & \leq & \Xi^*_{x,\gamma}(\alpha) \; \leq \; \alpha F^* + (1-\alpha) F(x_k). \QF
\ea
$$

Note that method (\ref{met-QN3}) can be endowed with a termination criterion, which supports its rate of convergence (\ref{eq-Quartic3}). Indeed, in this method, in parallel with the minimizing sequence $\{x_k\}_{k \geq 0}$, we can update the following lower bounds for the optimal value of problem (\ref{prob-QRMin}):
\beq\label{def-Crit}
\ba{rcl}
\xi^*_0 & = & \Xi^*_{x_0,\gamma}(1), \quad \xi^*_{k+1} \; = \; \max\Big\{ \xi^*_k, \; \Xi^*_{x_k,\gamma}(1) \Big\}, \quad k \geq 0.
\ea
\eeq
Clearly, under conditions of Theorem \ref{th-Quartic3},for all $k \geq 0$, we have
\beq\label{eq-XiGrow}
\ba{rcl}
(a): \; \xi^*_k & \leq & \xi^*_{k+1} \; \leq \; F^*, \quad (b): \; \xi_{k+1}^* \; \geq \; \Xi^*_{x_k,\gamma}(1).
\ea
\eeq
\BT\label{th-Crit}
Let the conditions of Theorem \ref{th-Quartic3} be satisfied. Then, for all $k \geq 0$, we have
\beq\label{eq-Crit}
\ba{rcl}
F(x_k) - F^* & \leq & F(x_k) - \xi_k^* \; \leq \; (1-\alpha)^k (F(x_0) - \xi^*_0).
\ea
\eeq
\ET
\proof
Indeed, by the same reasons as in the proof of Theorem \ref{th-Quartic2}, we can see that
$$
\ba{c}
F(x_{k+1}) - \xi^*_{k+1} \; \leq \; \alpha \, \Xi^*_{x_k,\gamma}(1) + (1-\alpha)F(x_k) - \xi^*_{k+1} \\
\\
\; \stackrel{(\ref{eq-XiGrow}{\rm b})}{\leq} \; (1-\alpha)(F(x_k) - \xi^*_{k+1})\; \stackrel{(\ref{eq-XiGrow}{\rm a})}{\leq} \; (1-\alpha)(F(x_k) - \xi^*_{k}). \QF
\ea
$$

Thus, if the inequality
\beq\label{eq-Check}
\ba{rcl}
F(x_k) - \xi_k^* & \leq & \epsilon
\ea
\eeq
is satisfied, then the point $x_k$ is an $\epsilon$-solution of problem (\ref{prob-QRMin}).

Finally, let us provide method (\ref{met-QN3}) with a reasonable value of parameter $\gamma$. In accordance to the relation (\ref{eq-Gamma3}), it could be found from the balance equation
$$
\ba{rcl}
\left( {3 \gamma - 1 \over 3 \gamma + 1} \right)^3 & = &  
{\mu \kappa(\gamma) \over L(1+2\gamma)} .
\ea
$$
However, it is quite complicated. Moreover, we need to emphasize the asymptotic dependence of the coefficient $\alpha$ of (\ref{eq-Gamma3}) in the condition number $q$. Therefore, we are going to use an approximate solution of the above equation.
\BT\label{th-Relax}
Let $q = {\mu \over L} \leq 1$. For coefficient
\beq\label{def-Gamma}
\ba{rcl}
\gamma \; = \; \gamma_{*} & \Def & {1 \over 3\left(1 - {3 \over 11}q^{1/3}\right)} \; \leq \; {11 \over 24},
\ea
\eeq
we have
\beq\label{eq-Order}
\ba{rcl}
{3 \gamma - 1 \over 3 \gamma + 1} & \leq &
\left[{ q \, \kappa(\gamma) \over 1+2\gamma}  \right]^{1/3}.
\ea
\eeq
\ET
\proof
Let $\gamma = \gamma_{*} \Def {1 \over 3(1 - \tau_* q^{1/3})}$ with $\tau_* = {3 \over 11}$. Denote $\lambda = {3 \gamma - 1 \over 3 \gamma} = \tau_*q^{1/3}$.
We need to prove that
$$
\ba{rcl}
(1+2\gamma)
\left({3 \gamma - 1 \over 3 \gamma + 1}\right)^3 & \leq &
q \left( (1-2 \gamma) - {16 \over 125}  \left(q^{-1}-1\right)\left({3 \gamma - 1 \over \gamma} \right)^3 \right)\\
\\
& = & q \left( (1-2 \gamma) + 16 \left({3 \over 5}\right)^3 \lambda^3 \right) - 16 \left({3 \over 5}\right)^3 \lambda^3.
\ea
$$
Since $\gamma \leq \half$, we have $\left({3 \gamma - 1 \over 3 \gamma + 1}\right)^3 = \left({3 \gamma \over 3 \gamma + 1}\right)^3 \lambda^3 \leq \left({3 \over 5} \right)^3 \lambda^3$. Thus, it is enough to prove that
$$
\ba{rcl}
\left[ (1 + 2 \gamma) \left({3 \over 5} \right)^3 + 16 \left({3 \over 5} \right)^3 \right] \lambda^3 & \leq & q (1 - 2 \gamma + 16 \lambda^3).
\ea
$$
However, $\lambda^3 = O(q)$. Hence, we can relax our goal up to the following target inequality:
\beq\label{eq-AuxRelax}
\ba{rcl}
18 \left({3 \over 5} \right)^3\lambda^3 & \leq & { q \over 12} \quad \Big( \refLE{def-Gamma} \; q (1 - 2 \gamma) \quad  \Big).
\ea
\eeq
By definition of $\lambda$, this is the inequality $18 \left({9 \over 55} \right)^3 \leq {1 \over 12}$, which is valid since $18 \cdot 12 = 6^3$.

Thus, inequality (\ref{eq-AuxRelax}) is valid and we prove the relation (\ref{eq-Order}).
\qed

Now, using $\gamma = \gamma_{*}$ in inequality (\ref{eq-QReg}), and denoting $\tau_* = {3 \over 11}$, we get the following lower bound:
\beq\label{eq-QLow}
\ba{rl}
& f(y) - f(x) - \la \nabla f(x), y - x \ra\\
\\
\geq &  {3 q^{1/3} \over 22} \la \nabla^2 f(x)(y-x), y - x \ra +{L (1+2 \gamma_*) \over 24} \| y - x \|^4  \left({ 3 \gamma_* - 1 \over 3 \gamma_* + 1} \right)^3\\
\\
\geq & {3 q^{1/3} \over 22} \la \nabla^2 f(x)(y-x), y - x \ra +{5 \mu \over 3 \cdot 24} \| y - x \|^4  \left({ \tau_* \over 2 - \tau_* q^{1/3}} \right)^3\\
\\
\geq & {3 q^{1/3} \over 22} \la \nabla^2 f(x)(y-x), y - x \ra +{5 \mu \over 3 \cdot 24} \left( {3 \over 22} \right)^3 \| y - x \|^4  .
\ea
\eeq

Relation (\ref{eq-Order}) can be also used for justifying the rate of convergence of method (\ref{met-QN3}]). Indeed, if we choose in (\ref{met-QN3}) $\gamma = \gamma_{*}$, then in accordance to Theorems \ref{th-Quartic3} and \ref{th-Relax}, the rate of convergence (\ref{eq-Quartic3}) is described by
$$
\ba{rcl}
\alpha & = & \alpha_{*} \; \Def \; {3 \gamma_{*} - 1\over 3 \gamma_{*} + 1}.
\ea
$$
Note that
$$
\ba{rcl}
\alpha_{*} & \refEQ{def-Gamma} & {3q^{1/3} \over 22 - 3 q^{1/3}}.
\ea
$$
Thus, in order to get $\epsilon$-solution of problem (\ref{prob-Quart}), the pure second-order method (\ref{met-QN3}) with $\gamma = \gamma_{*}$ needs
\beq\label{eq-Comp3}
\ba{c}
{22 \over 3 q^{1/3}} \ln {F(x_0) - F^* \over \epsilon}
\ea
\eeq
iterations at most. Note that this bound is approximately in four times worse than the bound (\ref{eq-CompRQN}) for efficiency of a specialized method (\ref{met-QN2}) as applied to the quartic polynomials.

\section{Applications in Convex Optimization}\label{sc-Conc}
\SetEQ

Let us show now that the notion of Q-regularity can be used for constructing new and efficient optimization schemes. In this section, we are working with Euclidean norm $\| x \| = \la B x, x \ra^{1/2}$ defined by a self-adjoint linear operator $B \succ 0$.

Consider first the problem (\ref{prob-QRMin}) with bounded fourth derivative of function $f(\cdot)$:
\beq\label{ass-Bound4}
\ba{rcl}
0 & \leq & D^4f(x)[h]^4 \; \leq \; L \| h \|^4, \quad \forall x \in S, \; h \in \E.
\ea
\eeq
Since $\mu = 0$, such a problem cannot be solved directly by method (\ref{met-QN3}). However, we can treat it by an appropriate regularization technique.

Indeed, let $\epsilon>0$ be our target accuracy in function value for an approximate solution of problem (\ref{prob-QRMin}). Suppose that we know an upper bound for the distance to the optimal solution:
\beq\label{eq-Dist}
\ba{rcl}
\| x_0 - x^* \| & \leq & R.
\ea
\eeq
Then we can define the following regularized version of the objective function:
\beq\label{def-Freg}
\ba{rcl}
F_{x_0,H}(x) & \Def & F(x) + {H \over 24} \| x - x_0 \|^4, \quad x \in \dom \psi,\\
\\
H & = &  {12 \, \epsilon \over R^4}.
\ea
\eeq
Note that this function is Q-regular with the following parameters:
\beq\label{eq-RConst}
\ba{rcl}
\mu_{\epsilon} & = & H, \quad L_{\epsilon} \; = \; L + H, \quad q_{\epsilon}^{-1} \; = \; 1 + {LR^4 \over 12 \epsilon}.
\ea
\eeq
Hence, since $F_{x_0,H}(x_0) = F(x_0)$ and $F^*_{x_0,H} \geq F^*$,
inequality (\ref{eq-Comp3}) provides us with the following bound on the number of iterations of method (\ref{met-QN3}) needed for finding an $(\epsilon/2)$-solution of the regularized problem:
\beq\label{eq-Sol2}
\ba{rcl}
{22 \over 3} \left[ 1 + {LR^4 \over 12 \epsilon} \right]^{1/3} \ln {2 (F(x_0) - F^*) \over \epsilon}.
\ea
\eeq
Let $\bar x$ be such a point, and $x^*$ be an optimal solution of problem (\ref{prob-QRMin}). Then
$$
\ba{rcl}
F(\bar x) & \leq & F_{x_0,H}(\bar x) \; \leq \; F^*_{x_0,H} + {\epsilon \over 2} \; \leq \; F_{x_0,H}(x^*) + {\epsilon \over 2} \\
\\
& = & F(x^*) + {\epsilon \over 2 R^4} \| x^* - x_0 \|^4 + {\epsilon \over 2} \; \refLE{eq-Dist} \; F^* + \epsilon.
\ea
$$
Thus, the point $\bar x$ is an $\epsilon$-solution of our initial problem (\ref{prob-QRMin}).

However, the notion of Q-regularity finds much more important applications in the framework of high-order proximal-point methods \cite{High,Masoud1}. In this approach, we solve the general problems of composite minimization in the form
\beq\label{prob-GMin}
\ba{rcl}
F^* & = & \min\limits_{x \in \dom \psi} \Big[ \; F(x) \; = \; f(x) + \psi(x) \; \Big],
\ea
\eeq
where $\psi(\cdot)$ is a simple closed convex function, and the smooth part $f(\cdot)$ of the objective function is convex and has bounded fourth derivative:
\beq\label{ass-Bound44}
\ba{rcl}
|\, D^4f(x)[h]^4 \,| & \leq & M_4 \| h \|^4, \quad \forall x \in (\dom \psi)_{1/2}, \; h \in \E.
\ea
\eeq

The problem (\ref{prob-GMin}) can be solved by different third-order proximal-point methods \cite{High,Prox2}, provided that we are able to compute in a reasonable time an approximate solution to the following auxiliary optimization problems:
\beq\label{prob-Prox}
\ba{c}
\min\limits_{x \in \dom \psi} \Big[\;  F_{\bar x, H}(x) = F(x) + {H \over 24} \| x - \bar x \|^4 \; \Big]
\ea
\eeq
with different prox-centers $\bar x \in \dom \psi$. Note that the objective function of this problem can be represented as follows:
$$
\ba{rcl}
F_{\bar x, H}(x) & = & f_{\bar x, H}(x) + \psi(x),\\
\\
f_{\bar x, H}(x) & = & f(x) + {H \over 24} \| x - \bar x \|^4.
\ea
$$
Hence, for $H$ big enough, the smooth part of the objective function in (\ref{prob-Prox}) is $Q$-regular. 
\BL\label{lm-Simple}
Let convex function $f(\cdot)$ satisfy condition (\ref{ass-Bound44}). Then for $H \geq M_4$, function $f_{\bar x,H}(\cdot)$ is Q-regular with parameters $\mu = H-M_4$ and $L = H+M_4$.
\EL
\proof
Indeed, since the norm $\| \cdot \|$ is Euclidean, for function $d(x) = {1 \over 24} \| x \|^4$ we have 
$$
\ba{rcl}
D^4d(x)[h]^4 & = & \| h \|^4, \quad \forall x, h \in \E.
\ea
$$ 
It remains to use inequality (\ref{ass-Bound44}).
\qed

Thus, for $H$ big enough, problem (\ref{prob-Prox}) can be efficiently solved by the Quartic Newton Method (\ref{met-QN3}). This method has linear rate of convergence in function value. Let us find an upper bound for the number of iterations, which is sufficient for finding an appropriate approximate solution of problem (\ref{prob-Prox}).

We need two auxiliary results. Let $\Phi(x) = \vf(x) + \psi(x)$, where $\psi(\cdot)$ is a simple closed convex function, and function $\vf(\cdot)$ is twice continuously differentiable on $\dom \psi$. Consider the following problem of composite minimization:
$$
\Phi_* \; = \; \min\limits_{x \in \dom \psi} \Phi(x).
$$ 
\BL\label{lm-Aux}
Let us fix a point $\bar x \in \dom \psi$ with $\Phi(\bar x) > \Phi_*$. Define the point
\beq\label{def-TPoint}
\ba{rcl}
T \; = \; T(\bar x)  & \Def & \arg\min\limits_{x \in \dom \psi} \Big[ \; \la \nabla \vf(\bar x), x \ra + {M \over 2} \| x - \bar x \|^2 + \psi(x) \; \Big],
\ea
\eeq
where $M \geq \hat M = \sup\limits_{x \in {\cal F}} \| \nabla^2 \vf(x) \|$ and ${\cal F} = \{ x \in \dom \psi: \; \Phi(x) \leq \Phi(\bar x) \}$. Then
\beq\label{eq-GNorm}
\ba{rcl}
\Phi'(T) & \Def & \nabla \vf(T) - \nabla \vf(\bar x) - M B(T-\bar x) \; \in \; \partial \Phi(T),\\
\\
(\| \Phi'(T) \|^*)^2 & \leq & 2M [\Phi(\bar x) - \Phi(T)].
\ea
\eeq
\EL
\proof
Consider the point
$\hat T = \arg\min\limits_{x \in {\cal F}} \Big[\; \la \nabla \vf(\bar x), x \ra + {M \over 2} \| x - \bar x \|^2 + \psi(x) \; \Big]$.
The first-order optimality condition for this problem is as follows:
\beq\label{eq-Opt1P}
\ba{rcl}
\la \nabla \vf(\bar x) + MB(\hat T - \bar x), x - \hat T \ra + \psi(x) & \geq & \psi(\hat T), \quad \forall x \in {\cal F}.
\ea
\eeq
If $\hat T = \bar x$, then for $x_* \in \mbox{Arg}\min\limits_{dom \psi} \Phi(x)$, we have
$$
\ba{rcl}
\Phi(x_*) & = & \vf(x_*) + \psi(x_*) \; \geq \; \vf(\bar x) + \la \nabla \vf(\bar x), x_* - \bar x \ra + \psi(x_*) \; \refGE{eq-Opt1P} \; \Phi(\bar x),
\ea
$$
and this contradicts to the assumptions of the theorem.

Thus, $\hat T \neq \bar x$.
Applying (\ref{eq-Opt1P}) to $x = \bar x$ and since both $\bar x$ and $\hat T$ belong to ${\cal F}$, we have
$$
\ba{rcl}
\Phi(\hat T) & = & \vf(\hat T) + \psi(\hat T)\\
\\
& \leq & \vf(\bar x) + \la \nabla \vf(\bar x), \hat T - \bar x \ra + {M \over 2} \| \hat T - \bar x \|^2 + \psi(\Hat T)\\
\\
& \refLE{eq-Opt1P} & \Phi(\bar x) - {M \over 2} \| \hat T - \bar x \|^2 \; < \; \Phi(\bar x).
\ea
$$
Hence, the functional boundary of the level set ${\cal F}$ is not active at $\hat T$ and therefore $\hat T = T$ as defined by (\ref{def-TPoint}).

Note that the first-order optimality condition for problem (\ref{def-TPoint}) is as follows:
\beq\label{eq-OptPT}
\ba{rcl}
\la \nabla \vf(\bar x) + MB(T - \bar x), x - T \ra + \psi(x) & \geq & \psi(T), \quad \forall x \in \dom \psi.
\ea
\eeq
This means that $\psi'(T) \Def - \nabla \vf(\bar x) - MB(T - \bar x) \in \partial \psi(T)$ and the inclusion 
$$
\ba{rcl}
\Phi'(T) & = & \nabla \vf(T) + \psi'(T) \in \partial \vf(T)
\ea
$$
is valid. On the other hand, since both $\bar x$ and $T$ belong to ${\cal F}$, we conclude that 
$$
\ba{rcl}
(\| \Phi'(T) \|^*)^2 & = & (\| \vf(T) - \nabla \vf(\bar x) - MB(T - \bar x) \|^*)^2\\
\\
& = & (\| \vf(T) - \nabla \vf(\bar x) \|^*)^2 - 2 M \la \vf(T) - \nabla \vf(\bar x),  T - \bar x \ra + M^2 \| T - \bar x \|^2\\
\\
& \refLE{eq-GDiff} & M^2 \| T - \bar x \|^2.
\ea
$$
This gives us the bound (\ref{eq-GNorm}).
\qed

We are going to apply Lemma \ref{lm-Aux} in the situation when function $\vf(\cdot)$ in Q-regular. Hence, we need to use an explicit bound for its second derivative.
\BL\label{lm-HBound}
Let function $\vf(\cdot)$ be Q-regular on $\dom \psi$ with constants $L \geq \mu > 0$. Then, for any $x \in {\cal F} = \{ x \in \dom \psi: \; \Phi(x) \leq \Phi(\bar x) \}$ we have
\beq\label{eq-HatM}
\ba{rcl}
\| \nabla^2 \vf(x) \| & \leq & \hat M \Def 4 \| \nabla^2 \vf(\bar x) \| + {2 \over 3} LD^2,
\ea
\eeq
where $D \Def 2 \left[ {72 \over \mu} (\Phi(\bar x) - \Phi_*)\right]^{1/4}$.
\EL
\proof
Since $\mu > 0$, function $\Phi(\cdot)$ is uniformly convex. Let us fix $\gamma= {1 \over 3}$.
By inequality (\ref{eq-QReg}), we get
$$
\ba{rcl}
{\mu \over 72} \| \bar x - x_* \|^4 & \leq & \Phi(\bar x) - \Phi(x_*)  \; \Def \; \Delta,
\ea
$$
where $x_* = \arg\min\limits_{x \in \dom \psi} \Phi(x)$. Thus, $\| \bar x - x_* \| \leq R \Def \left[ {72 \over \mu} \Delta \right]^{1/4}$. Note that the same bound on the distance to the optimum is valid for any point from the level set ${\cal F}$.
Hence, for any $y \in {\cal F}$ we have $\| y - \bar x \| \leq D = 2 R$.

In view of (\ref{def-QRegU}), for any $y \in \dom \psi$ we have
$$
\ba{rcl}
0 & \preceq & \nabla^2 \vf(\bar x - \gamma(y-\bar x)) \; \preceq \; \nabla^2 \vf(\bar x) - \gamma D^3\vf(\bar x)[y - \bar x] + {\gamma^2 \over 2} L \| y - \bar x \|^2 B.
\ea
$$
Therefore, for any $y \in {\cal F}$, we get the following bound:
$$
\ba{rcl}
\nabla^2 \vf(y) & \preceq & \nabla^2 \vf(\bar x) + D^3 \vf(\bar x)[y-\bar x] + \half L \| y - \bar x \|^2 B\\
\\
& \preceq & {1 + \gamma \over \gamma} \nabla^2 \vf(\bar x) + {1 + \gamma \over 2} L \| y - \bar x \|^2 B.
\ea
$$
Thus, $\hat M = {1 + \gamma \over \gamma} \| \nabla^2 \vf(\bar x) \| + {1 + \gamma \over 2} LD^2$ is an upper bound for the norms of Hessians of function $\vf(\cdot)$ on the set ${\cal F}$.
\qed

Let us show now how we can use all this machinery in the framework of high-order proximal-point methods of degree three. Recall that an approximation $\hat x \in \dom \psi$ of the exact minimum of problem (\ref{prob-Prox}) is acceptable if it satisfies the following condition:
\beq\label{eq-Accept}
\ba{rcl}
\exists \hat g \in \partial \psi(\hat x): \quad \| \nabla F_{\bar x, H}(\hat x) + \hat g \|^* & \leq & \beta \| \nabla f(\hat x) + \hat g \|^*,
\ea
\eeq
where $\beta \in (0,1)$ is the tolerance parameter (see \cite{Masoud1}; compare with \cite{High}). Note that at the exact solution of problem (\ref{prob-Prox}) we can enforce the left-hand side of this condition to be zero. 

Let us estimate complexity of finding such a point by QRNM (\ref{met-QN3}). The whole process consists of two stages, dependent on positive parameteres $\delta$ and $M$.
\BI
\item{\bf Stage 1.}
Suppose that we know the constant $M_4$ in (\ref{ass-Bound44}). Then, taking $H = 2M_4$, we can ensure Q-regularity of function $f_{\bar x,h}(\cdot)$ with parameters
$$
\ba{rcl}
\mu & = & M_4, \quad L \; = \; 3 M_4, \quad q \; = \; {1 \over 3}
\ea
$$
(see Lemma \ref{lm-Simple}). Hence, choosing in (\ref{met-QN3}) $\gamma \refEQ{def-Gamma} \gamma_*$, we get the minimizing sequence convergent with the following rate:
\beq\label{eq-RateAux}
\ba{rcl}
F_{\bar x, H}(x_k) - \xi_k^* & \refLE{eq-Crit} & \left( 1 - \alpha_* \right)^k (F(\bar x) - F^*),
\ea
\eeq
where $\xi_k^*$ are computable lower bounds for the optimal value of (\ref{prob-Prox}) updated in accordance to (\ref{def-Crit}). This stage is terminated when $F_{\bar x, H}(x_k) - \xi_k^* \leq \delta$.
\item{\bf Stage 2.} Compute point $\hat x = T(x_k)$ using the iteration (\ref{def-TPoint}) with $\vf(\cdot) = f_{\bar x, H}(\cdot)$ and an appropriate constant $M$.
\EI

Let us point out the values of parameters $\delta$ and $M$ ensuring validity of the condition~(\ref{eq-Accept}). In view of Lemma~\ref{lm-HBound}, we can take
$$
\ba{rcl}
M & = & 4 \| \nabla^2 f(\bar x) \| + 2 D^2, \quad D \; = \; 2 \left[ {72 \over M_4} (F(\bar x) - F^*)\right]^{1/4}.
\ea
$$
Thus, we come to the following conclusion.
\BL\label{lm-Bound}
Let problem (\ref{prob-GMin}) satisfy the following assumptions:
\beq\label{ass-GMin}
\ba{rl}
(a): & \| \nabla^2 f(x) \| \; \leq \; M_2, \quad x \in \dom \psi, \\
\\
(b): & F(x) - F(y) \; \leq \; \Delta_F, \quad x, y \in \dom \psi.
\ea
\eeq
Then we can take 
\beq\label{def-M}
\ba{rcl}
M & = & 4 M_2 + 48 \left[ {2  \Delta_F \over M_4} \right]^{1/2}.
\ea
\eeq
\EL

Suppose now that $F(\hat x) - F^* \geq \epsilon$. Then in the condition (\ref{eq-Accept}) we have 
$$
\ba{rcl}
 \| \nabla f(\hat x) + \hat g \|^* & \geq & (F(\hat x) - F^*)/ D_{\psi} \; \geq \; \epsilon / D_{\psi},
\ea
$$
where $D_{\psi} < + \infty$ is the diameter of $\dom \psi$. Hence, in view of inequality (\ref{eq-GNorm}), in order to satisfy condition (\ref{eq-Accept}), it is sufficient to ensure
$$
\ba{rcl}
2 M (F_{\bar x, H}(x_k) - \xi_k^* ) & \leq & \left({\beta \, \epsilon \over D_{\psi}} \right)^2.
\ea
$$
Taking into account the rate of convergence (\ref{eq-RateAux}), we get the following upper bound for the number of iterations of method (\ref{met-QN3}) for solving the auxiliary problem (\ref{prob-Prox}):
\beq\label{eq-TotIter}
\ba{c}
{22 \over 3^{2/3}} \ln {2 M \Delta_F D_{\psi}^2 \over \beta^2 \epsilon^2}
\ea
\eeq
Having this machinery at hand, we can design different second-order methods for solving the problem (\ref{prob-GMin}) with the rate of convergence $O(k^{-p})$ with $p = 3, 4, 5$ (see \cite{Masoud1,Masoud2}).

\end{document}